\documentclass{amsart}

\newtheorem{theorem}{Theorem}[section]
\newtheorem{lemma}[theorem]{Lemma}

\usepackage[dvips]{graphics}

\theoremstyle{definition}

\theoremstyle{remark}

\theoremstyle{conjecture}

\theoremstyle{corollary}
\newtheorem{corollary}[theorem]{Corollary}

\numberwithin{equation}{section}

\begin{document}

\title{Matchings and entropies of cylinders}

\author{Weigen Yan}
\address{School of Sciences, Jimei University,
Xiamen 361021, China}
\curraddr{Institute of Mathematics, Academia Sinica, Taipei 11529,
Taiwan} \email{wgyan@math.sinica.edu.tw;\ \ weigenyan@263.net}
\thanks{The first author was supported in part by FMSTF Grant \#2004J024 and by NSFF Grant \#E0540007.}

\author{Yeong-Nan Yeh}
\address{Institute of Mathematics, Academia Sinica, Taipei 11529,
Taiwan} \email{mayeh@math.sinica.edu.tw}
\thanks{The second author was supported in part by NSC Grant \#95-2115-M-001-009.}

\author{Fuji Zhang}
\address{Corresponding author, School of
Mathematical Science, Xiamen University, Xiamen 361005, China}
\email{fjzhang@jingxian.xmu.edu.cn}
\thanks{The third author was supported in part by NSFC Grant \#10671162.}

%    General info
%; Secondary 46E25, 20C20}

%\date{January 1, 2001 and, in revised form, June 22, 2001.}



\begin{abstract}
The enumeration of perfect matchings of graphs is equivalent to
the dimer problem which has applications in statistical physics. A
graph $G$ is said to be $n$-rotation symmetric if the cyclic group
of order $n$ is a subgroup of the automorphism group of $G$.
Jockusch (Perfect matchings and perfect squares, J. Combin. Theory
Ser. A, 67(1994), 100-115) and Kuperberg (An exploration of the
permanent-determinant method, Electron. J. Combin., 5(1998), \#46)
proved independently that if $G$ is a plane bipartite graph of
order $N$ with $2n$-rotation symmetry, then the number of perfect
matchings of $G$ can be expressed as the product of $n$
determinants of order $N/2n$. In this paper we give this result a
new presentation. We use this result to compute the entropy of a
bulk plane bipartite lattice with $2n$-notation symmetry. We
obtain explicit expressions for the numbers of perfect matchings
and entropies for two types of cylinders. Using the results on the
entropy of the torus obtained by Kenyon, Okounkov, and Sheffield
(Dimers and amoebae, Ann. Math. 163(2006), 1019--1056) and by
Salinas and Nagle (Theory of the phase transition in the layered
hydrogen-bonded $SnCl^2\cdot 2H_2O$ crystal, Phys. Rev. B,
9(1974), 4920--4931), we show that each of the cylinders
considered and its corresponding torus have the same entropy.
Finally, we pose some problems.
\end{abstract}

\maketitle

%\section*{This is an unnumbered first-level section head}
%This is an example of an unnumbered first-level heading.

%% The correct journal style for \specialsection is all uppercase; a known bug
%% in amsart.cls prevents this, so input must be uppercase until it is fixed.
%\specialsection*{This is a Special Section Head}
%\specialsection*{THIS IS A SPECIAL SECTION HEAD}
%This is an example of a special section head%
%%%%%%%%%%%%%%%%%%%%%%%%%%%%%%%%%%%%%%%%%%%%%%%%%%%%%%%%%%%%%%%%%%%%%%%%
%\footnote{Here is an example of a footnote. Notice that this
%footnote text is running on so that it can stand as an example of
%how a footnote with separate paragraphs should be written.
%\par
%And here is the beginning of the second paragraph.}%
%%%%%%%%%%%%%%%%%%%%%%%%%%%%%%%%%%%%%%%%%%%%%%%%%%%%%%%%%%%%%%%%%%%%%%%%

\section{Introduction}
%This is an example of a numbered first-level heading.

%\subsection{This is a numbered second-level section head}
%This is an example of a numbered second-level heading.

%\subsection*{This is an unnumbered second-level section head}
%This is an example of an unnumbered second-level heading.

%\subsubsection{This is a numbered third-level section head}
%This is an example of a numbered third-level heading.

%\subsubsection*{This is an unnumbered third-level section head}
%This is an example of an unnumbered third-level heading.
An automorphism of a graph $G$ is a graph isomorphism with itself,
i.e., a bijection from the vertex set of $G$ to itself such that
the adjacency relation is preserved. The set of automorphisms
defines a permutation group on the vertices of $G$, which is
called the automorphism group of the graph $G$. The automorphism
group of a graph $G$ characterizes its symmetries, and is
therefore very useful for simplifying the computation of some of
its invariants. One of the most successful results in this
direction is to use the character of the automorphism group of a
graph $G$ to compute its characteristic polynomial (see for
example Cvetkovi\'c {\it et al.} \cite{CDS80}). Recall that a
cyclic group is a group that can be generated by a single element
$x$ (the group generator). A cyclic group of order $n$ is denoted
by $\mathcal C_n$, and its generator $x$ satisfies $x^n=1$, where
$1$ is the identity element. A graph $G$ is said to be
$n$-rotation symmetric if the cyclic group $\mathcal C_n$ is a
subgroup of the automorphism group of $G$. We call $\mathcal C_n$
the rotational subgroup of an $n$-rotation symmetric graph $G$. We
say that a plane graph $G$ is reflective symmetric if it is
invariant under the reflection across some straight line.

A perfect matching of a graph $G$ is a collection of
vertex-disjoint edges that are collectively incident to all
vertices of $G$. Let $M(G)$ denote the number of perfect matchings
of a simple graph $G$.  Problems involving enumeration of perfect
matchings have been examined extensively not only by
mathematicians (see for example
\cite{Ciuc97,CEP96,CKP01,Fisher61,KOS06AM,Prop99,Stem90}) but also
by physicists and chemists (see for example
\cite{GC89,Kast61,Kast63,Kast67,TF61,WFY68}). In 1961, Kasteleyn
\cite{Kast61} found a formula for the number of perfect matchings
of an $m\times n$ quadratic lattice graph. Temperley and Fisher
\cite{TF61} used a different method and arrived at the same result
at almost exactly the same time. Both lines of calculation showed
that the logarithm of the number of perfect matchings, divided by
$\frac{mn}{2}$, converges to $2c/\pi\approx 0.5831$ as
$m,n\rightarrow \infty$, where $c$ is Catalan's constant. This
limit is called the entropy of the quadratic lattice graph and the
corresponding problem was called the dimer problem by the
statistical physicists. In 1992, Elkies {\it et al.} \cite{EKLP92}
studied the enumeration of perfect matchings of regions called
Aztec diamonds, and showed that
%that the logarithm of the number of
%perfect matchings, divided by the number of edges in each perfect
%matching, converges to the smaller number
the entropy equals $\frac{\log 2}{2}\approx 0.35$. Cohn, Kenyon,
and Propp \cite{CKP01} demonstrated that the behavior of random
perfect matchings of large regions $R$ was determined by a
variational (or entropy maximization) principle, as was
conjectured in Section 8 of \cite{EKLP92}, and they gave an exact
formula for the entropy of simply-connected regions of arbitrary
shape. Particularly, they showed that computation of the entropy
is intimately linked with an understanding of long-range
variations in the local statistics of random domino tilings.
Kenyon, Okounkov, and Sheffield \cite{KOS06AM} considered the
problem enumerating perfect matchings of the doubly-period
bipartite graph on a torus, which generalized the results in
\cite{CKP01}. They proved that the number of perfect matchings of
the doubly-period plane bipartite graph $G$ can be expressed in
terms of four determinants and they expressed the entropy of $G$
as a double integral.

The dimer problem on hexagonal lattice graphs has also been
examined in the past \cite{Elser84,WFY68,WWB89}. It can be seen as
equivalent to the combinatorial problem of ``plane partitions"
(see MacMahon \cite{MacM}). The $a, b, c$ semiregular hexagon is
the hexagon whose side lengths are in cyclic order, $a, b, c, a,
b, c$. Lozenge tilings of this region are in correspondence with
plane partitions with at most $a$ rows, at most $b$ columns, and
no part exceeding $c$. MacMahon \cite{MacM} showed that the number
of such plane partitions ({\it i.e.}, the number of perfect
matchings of the corresponding hexagonal lattice graph $G(a,b,c)$)
equals
$$\prod_{i=1}^a\prod_{j=1}^b\prod_{k=1}^c\frac{i+k+k-1}{i+j+k-2}.$$
Elser \cite{Elser84} obtained the entropy of $G(a,a,a)$ which is
approximately 0.2616. The dimer problem on some other types of
plane lattices are also considered (see Wu \cite{WFY03} and the
references cited therein).

Jockusch \cite{Jock94} first used a combinatorial method to
obtained a product formula for the number of perfect matchings of
a type of rotational symmetric bipartite graphs. Kuperberg
\cite{Kupe98} used the representation theory of groups to obtain
independently a product formula which can be equivalent to the
Jockusch's one. Ciucu \cite{Ciuc97}, and Yan and Zhang
\cite{YZ04,YZ06,ZY03} considered the enumeration of perfect
matchings of general graphs (lattices) with a certain type of
reflective symmetry. Ciucu \cite{Ciuc97} obtained a basic
factorization theorem for the number of perfect matchings of plane
bipartite graphs with reflective symmetry. Yan and Zhang
\cite{YZ04,YZ06} extended Ciucu's result to general plane graphs
(not necessary to be bipartite) with reflective symmetry by using
Pfaffians.

The current paper deals with the enumerating problem for perfect
matchings of cylinders which can be regarded as the graphs with
rotational symmetry. In the next section, we introduce the
Pfaffian method for enumerating perfect matchings of plane graphs
(see for example \cite{Fisher61,Kast61,Stem90,TF61}). In Section
3, we give a new presentation of the product formula for the
number of perfect matchings of plane bipartite graphs with
$2n$-rotation symmetry, which was found independently by Jockusch
\cite{Jock94} and by Kuperberg \cite{Kupe98}. We use this result
to compute the entropy of a bulk plane bipartite lattice with
$2n$-notation symmetry. In Section 4, we obtain the explicit
expressions for the numbers of perfect matchings and entropies for
two types of cylinders (bipartite graphs with the
cylinder-boundary condition). Using the theorem obtained by
Kenyon, Okounkov, and Sheffield \cite{KOS06AM}, we compute the
entropy of a type of toruses (bipartite graphs with the
torus-boundary condition). On the other hand, the entropy of
another type of toruses was computed by Salinas and Nagle
\cite{SN74}. Our results show that each of the cylinders
considered and its corresponding torus have the same entropy (a
similar result for the plane quadratic lattices had been obtained
by Kasteleyn \cite{Kast61}). Finally, Section 5 poses two open
problems.
%As applications, in the last section
%we obtain the explicit expressions of the numbers of perfect
%matchings and entropies per dimer for two types of bulk
%tessellation cylinders---the 4-6-8 cylinder and the 4-8 cylinder.
%%%%%%%%%%%%%%%%%%%%%%%%%%%%%%%%%%%%%%%%%%%%%%%%%%%%%%%%%%%%%%%%%%%%%%%%

\section{Pfaffians}

Let $B=(b_{ij})_{2n\times 2n}$ be a skew symmetric matrix of order
$2n$. For each partition
$P=\{\{i_1,j_1\},\{i_2,j_2\},\ldots,\{i_n,j_n\}\}$ of the set
$\{1,2,\ldots,2n\}$ into pairs, form the expression
$$b_P=sgn(i_1j_1i_2j_2\ldots i_nj_n)b_{i_1j_1}b_{i_2j_2}\ldots b_{i_nj_n},$$
where $sgn(i_1j_1i_2j_2\ldots i_nj_n)$ denotes the sign of the
permutation $i_1j_1i_2j_2\ldots i_nj_n$. Note that $b_P$ depends
neither on the order in which the classes of the partition are
listed nor on the order of the two elements of a class. So $b_P$
indeed depends only on the choice of the partition $P$. The
Pfaffian of the skew matrix $B$ (see \cite{LP86}), denoted by
$Pf(B)$, is defined as
$$Pf(B)=\sum_{P}b_P,$$
where the summation ranges over all partitions of
$\{1,2,\ldots,2n\}$ which have the form of $P$.
%%%%%%%%%%%%%%%%%%%%%%%%%%%%%%%%%%%%%%%%%%
%%%%%%%%%%%%% Theorem 2.1
\begin{theorem} [The Cayley's Theorem, \cite{LP86}]
Let $B=(b_{ij})_{2n\times 2n}$ be a skew symmetric matrix of order
of $2n$. Then $$\det(B)=Pf(B)^2.$$
\end{theorem}

The Pfaffian method for enumerating perfect matchings of plane
graphs was independently discovered by Fisher \cite{Fisher61},
Kasteleyn \cite{Kast61}, and Temperley \cite{TF61}. Given a plane
graph $G$, the method produces a matrix $A$ such that the number
of perfect matchings of $G$ can be expressed by the determinant of
matrix $A$. By using this method,  Fisher \cite{Fisher61},
Kasteleyn \cite{Kast61}, and Temperley \cite{TF61} solved
independently a famous problem on enumerating perfect matchings of
an $m\times n$ quadratic lattice graph in statistical
physics--Dimer problem. Given a simple graph $G=(V(G), E(G))$ with
vertex set $V(G)=\{v_1,v_2,\ldots,v_N\}$, let $G^e$ be an
arbitrary orientation. The skew adjacency
matrix of $G^e$, denoted by $A(G^e)$, is defined as follows:\\
$$
A(G^e)=(a_{ij})_{n\times n},$$ where
$$
a_{ij}=\left\{\begin{array}{rl}
1 & \ \mbox{if}\ \  (v_i,v_j)\ \mbox{is an arc of}\ G^e,\\
-1 & \ \mbox{if} \ \  (v_j,v_i)\ \mbox{is an arc of}\ G^e,\\
0 & \ \mbox{otherwise}.
\end{array}
\right.$$  Obvious, $A(G^e)$ is a skew symmetric matrix.

If $D$ is an orientation of a graph $G$ and $C$ is a cycle of even
length, we say that $C$ is oddly oriented in $D$ if $C$ contains
odd number of edges that are directed in $D$ in the direction of
each orientation of $C$. We say that $D$ is a Pfaffian orientation
of $G$ if every nice cycle of even length of $G$ is oddly oriented
in $D$ (a cycle $C$ of $G$ is called be nice if the induced
subgraph $G-C$ of $G$ by $V(G)\backslash V(C)$ has perfect
matchings). It is well known that if a graph $G$ contains no
subdivision of $K_{3,3}$ then $G$ has a Pfaffian orientation (see
Little \cite{Litt75}). McCuaig \cite{McCu04}, and McCuaig,
Robertson et al \cite{MRST97}, and Robertson, Seymour et al
\cite{RST99} found a polynomial-time algorithm to show whether a
bipartite graph has a Pfaffian orientation. Stembridge
\cite{Stem90} proved that the number (or generating function) of
nonintersecting $r$-tuples of paths from a set of $r$ vertices to
a specified region in an acyclic digraph $D$ can, under favorable
circumstances, be expressed as a Pfaffian. See a survey of
Pfafffian orientations of graphs in Thomas \cite{Thomas06}.
%Some related work, see
%for example papers [13$-$16]

%%%%%%%%%%%%%%%%%%%%%%%%%%%%%%%%%%%%%%%%%%
%%%%%%%%%%%%% Lemma 2.2
\begin{lemma}[Lov\'asz et al. \cite{LP86}]
If $G^e$ is a Pfaffian orientation of a graph $G$, then
$$
M(G)=|Pf(A(G^e))|,
$$
where $A(G^e)$ is the skew adjacency matrix of $G^e$.
\end{lemma}

%%%%%%%%%%%%%%%%%%%%%%%%%%%%%%%%%%%%%%%%%%
%%%%%%%%%%%%% Lemma 2.3
\begin{lemma}[Lov\'asz {\it et al.} \cite{LP86}]
If $G$ is a plane graph and $G^e$, an orientation of $G$ such that
every boundary face--except possibly the infinite face--has an odd
number of edges oriented clockwise, then in every cycle the number
of edges oriented clockwise is of opposite parity to number of
vertices of $G^e$ inside the cycle. Consequently, $G^e$ is a
Pfaffian orientation.
\end{lemma}

%%%%%%%%%%%%%%%%%%%%%%%%%%%%%%%%%%%%%%%%%%
%%%%%%%%%%%%% Lemma
\begin{lemma}[\cite{Fisher61,Kast61,TF61}]
Every plane graph has a Pfaffian orientation satisfying the
condition in Lemma 2.3.
\end{lemma}
%%%%%%%%%%%%%%%%%%%%%%%%%%%%%%%%%%%%%%%%%%
%%%%%%%%%%%%% Figure 1
%%%%%%%%%%%%%%%%%%%%%%%%%%%%%%%%%%%%%%%%%%
\begin{figure}[htbp]
  \centering
 \scalebox{0.7}{\includegraphics{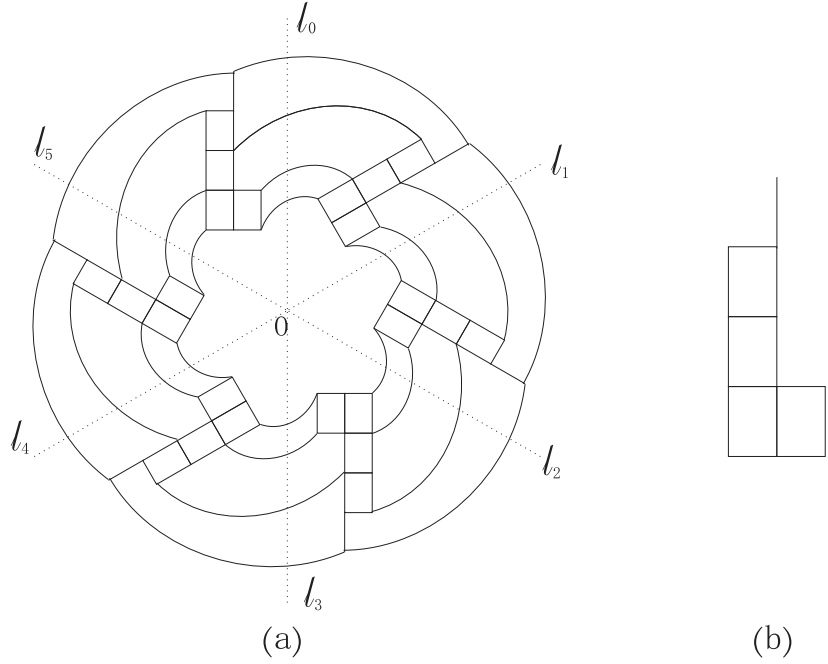}}
  \caption{\ (a)\ A connected plane bipartite graph $G$
with $6$-rotation symmetry, where the local width $w(G)$ of $G$
equals $2$.\ (b)\ One component $G_0$ of $G$.}
\end{figure}

%%%%%%%%%%%%%%%%%%%%%%%%%%%%%%%%%%%%%%%%%%%%%%%%%%%%%%%%%%%%%%%%%%%%%%%%

\section{A PRODUCT THEOREM}

Let $G=(V(G),E(G))$ be a simple connected plane bipartite graph of
order $N$ with $2n$-rotation symmetry, which has symmetric axes
$\ell_0$, $\ell_1$, $\ldots$, and $\ell_{2n-1}$ passing through
the rotation center $O$. Figure 1(a) shows an example of a
$6$-rotation symmetric graph. If $G$ has some vertices lying on
its symmetric axes $\ell_i$'s, then we can rearrange vertices of
$G$ such that there are no vertex lying on its new symmetric axes.
Hence we can assume that there exists no vertex of $G$ lying on
$\ell_0$, $\ell_1$, $\ldots$, $\ell_{2n-1}$. If we delete all the
edges intersected by $\ell_0$, $\ell_1$, $\ldots$, and
$\ell_{2n-1}$, then $2n$ isomorphic components
$G_i=(V(G_i),E(G_i))$ for $0\leq i\leq 2n-1$ are obtained, where
$V(G_i)=\{v_1^{(i)},v_2^{(i)},\ldots,v_{N/2n}^{(i)}\}$ is the
vertex set of $G_i$. Note that $G$ is a plane graph. Hence there
exists no edges between $G_i$ and $G_j$ for $j\neq i-1, i+1$ (mod
$2n$). That is, we can arrange all $G_i$'s on a circle clockwise
(see Figure 1(a)). Without loss of generality, we assume that the
set of edges of $G$ between $G_i$ and $G_{i+1}$ is exactly
$\{v_{r_k}^{(i)}v_{s_k}^{(i+1)}|k=1,2,\ldots,p\}$ (see Figure
1(a)), where $G_{2n}=G_0$, $v_{s_k}^{(2n)}=v_{s_{k}}^{(0)}$, and
$v_{r_k}^{(2n)}=v_{r_{k}}^{(0)}$, and the generator $x$ of the
rotation subgroup $\mathcal C_{2n}$ of the automorphism group of
$G$ maps $v_{r_k}^{(i)}$ (resp. $v_{s_k}^{(i)}$) to
$v_{r_k}^{(i+1)}$ (resp. $v_{s_k}^{(i+1)}$) for $1\leq k \leq p$,
respectively. The local width of a $2n$-rotation symmetric graph
$G$ with $2n$ isomorphic components $G_i$'s, denoted $w(G)$, is
defined to be the number of edges in $G_i$ lying on the boundary
face of $G$ which contains the rotation center $O$. For the graph
$G$ illustrated in Figure 1(a), $w(G)=2$. For the underlying graph
$G$ (which is a $6$-rotation symmetric graph) of the digraph
illustrated in Figure 4(a), $w(G)=1$.

Note that $G_0$ is a plane graph. By Lemma 2.4, there exists a
Pfaffian orientation $G_0^e$ of $G_0$ satisfying the condition in
Lemma 2.3. For the corresponding graph $G_0$ illustrated in Figure
1(b), the orientation $G_0^e$ shown in Figure 2(a) is a Pfaffian
orientation satisfying the condition in Lemma 2.3.
%%%%%%%%%%%%%%%%%%%%%%%%%%%%%%%%%%%%%%%%%%
%%%%%%%%%%%%% Figure 2
%%%%%%%%%%%%%%%%%%%%%%%%%%%%%%%%%%%%%%%%%%
\begin{figure}[htbp]
  \centering
 \scalebox{0.9}{\includegraphics{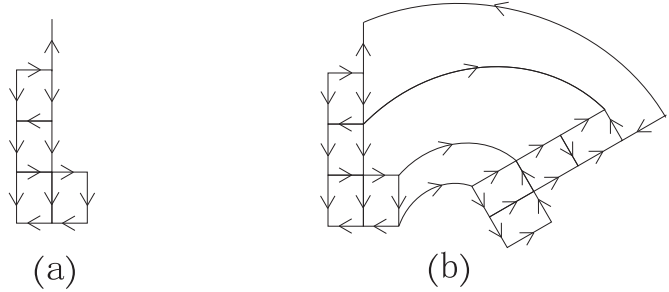}}
  \caption{\ (a)\ A Pfaffian orientation $G_0^e$ of $G_0$.\ (b)\ The corresponding Pfaffian
  orientation $G_{0,1}^e$ of $G_{0,1}$.}
\end{figure}
Now we define a Pfaffian orientation $G^e$ from $G_0^e$ as
follows. If we reverse the orientation of each arc of $G_0^e$,
then we can obtain an orientation of $G_0$ satisfying conditions
in Lemma 2.3, denoted by $G_0^{-e}$. Hence, $G_0^{-e}$ is the
converse of $G_0^{e}$. Note that $G_0$ is bipartite. Therefore,
$G_0^{-e}$ is a Pfaffian orientation of $G_0$ satisfying the
condition in Lemma 2.3. Since all $G_i$'s are isomorphic, we
define the orientations of all $G_{2i}$'s to be the same as
$G_0^e$ and the orientations of all $G_{2i+1}$'s to be same as
$G_0^{-e}$ for $0\leq i\leq n-1$, respectively. In order to obtain
an orientation $G^e$ of $G$, we need to orient all edges with the
form of $v_{r_k}^{(i)}v_{s_k}^{(i+1)}$ in $G$. Note that the
subgraph $G_{0,1}$ of $G$ induced by $V(G_0)\cup V(G_1)$ is a
plane graph. Define the orientations of $G_0$ and $G_1$ in
$G_{0,1}$ to be $G_0^e$ and $G_0^{-e}$, respectively, and orient
edges between $G_0$ and $G_1$ one by one such that every face
between $G_0$ and $G_1$ has an odd number of edges oriented
clockwise. Hence we have obtained an orientation $G_{0,1}^e$ of
$G_{0,1}$ satisfying the condition in Lemma 2.3 such that the the
induced orientations of $G_0$ and $G_1$ in $G_{0,1}^e$ are exactly
$G_0^e$ and $G_0^{-e}$. For the graph $G$ illustrated in Figure
1(a) and the orientation $G_0^e$ shown in Figure 2(a), the
corresponding orientation $G_{0,1}^e$ is pictured in Figure 2(b).

Hence we have oriented the edges $v_{r_k}^{(0)}v_{s_k}^{(1)}$ for
$1\leq k\leq p$ in $G$. Define the orientations of edges
$v_{r_k}^{(i)}v_{s_{k}}^{(i+1)}$ for $1\leq i\leq 2n-2$ to be from
$v_{r_k}^{(i)}$ to $v_{s_{k}}^{(i+1)}$ if the orientation of
$v_{r_k}^{(0)}v_{s_k}^{(1)}$ is from $v_{r_k}^{(0)}$ to
$v_{s_{k}}^{(1)}$, and the orientations of edges
$v_{r_k}^{(i)}v_{s_{k}}^{(i+1)}$ to be from $v_{s_k}^{(i+1)}$ to
$v_{r_{k}}^{(i)}$ otherwise. So far, we have oriented all edges of
$G$ except the edges $v_{r_k}^{(2n-1)}v_{s_{k}}^{(0)}$ for $1\leq
k\leq p$. In order to orient the edges
$v_{r_k}^{(2n-1)}v_{s_{k}}^{(0)}$ such that $G^e$ is a Pfaffian
orientation, we must distinguish the following four cases:

{\it Case (i)}.\ \ both $n$ and $w(G)$ (the local width of $G$)
are even;

{\it Case (ii)}.\ \ $n$ is odd and $w(G)$ is even;

{\it Case (iii)}.\ \ $n$ is even and $w(G)$ is odd;

{\it Case (iv)}.\ \ both $n$ and $w(G)$ are odd.

For {\it Cases (i)}, {\it (ii)}, and {\it (iii)}, define the
orientations of edges $v_{r_k}^{(2n-1)}v_{s_{k}}^{(0)}$ to be from
$v_{r_k}^{(2n-1)}$ to $v_{s_{k}}^{(0)}$ if the orientation of
$v_{r_k}^{(0)}v_{s_k}^{(1)}$ is from $v_{s_k}^{(1)}$ to
$v_{r_{k}}^{(0)}$, and the orientations of edges
$v_{r_k}^{(2n-1)}v_{s_{k}}^{(0)}$ to be from $v_{s_k}^{(0)}$ to
$v_{r_{k}}^{(2n-1)}$ otherwise (see Figure 3(a)).

For {\it Case (iv)}, define the orientations of edges
$v_{r_k}^{(2n-1)}v_{s_{k}}^{(0)}$ to be from $v_{r_k}^{(2n-1)}$ to
$v_{s_{k}}^{(0)}$ if the orientation of
$v_{r_k}^{(0)}v_{s_k}^{(1)}$ is from $v_{r_{k}}^{(0)}$ to
$v_{s_k}^{(1)}$, and the orientations of edges
$v_{r_k}^{(2n-1)}v_{s_{k}}^{(0)}$ to be from $v_{s_k}^{(0)}$ to
$v_{r_{k}}^{(2n-1)}$ otherwise (see Figure 4(a)).
%%%%%%%%%%%%%%%%%%%%%%%%%%%%%%%%%%%%%%%%%%
%%%%%%%%%%%%% Figure 3
%%%%%%%%%%%%%%%%%%%%%%%%%%%%%%%%%%%%%%%%%%
\begin{figure}[htbp]
  \centering
 \scalebox{0.7}{\includegraphics{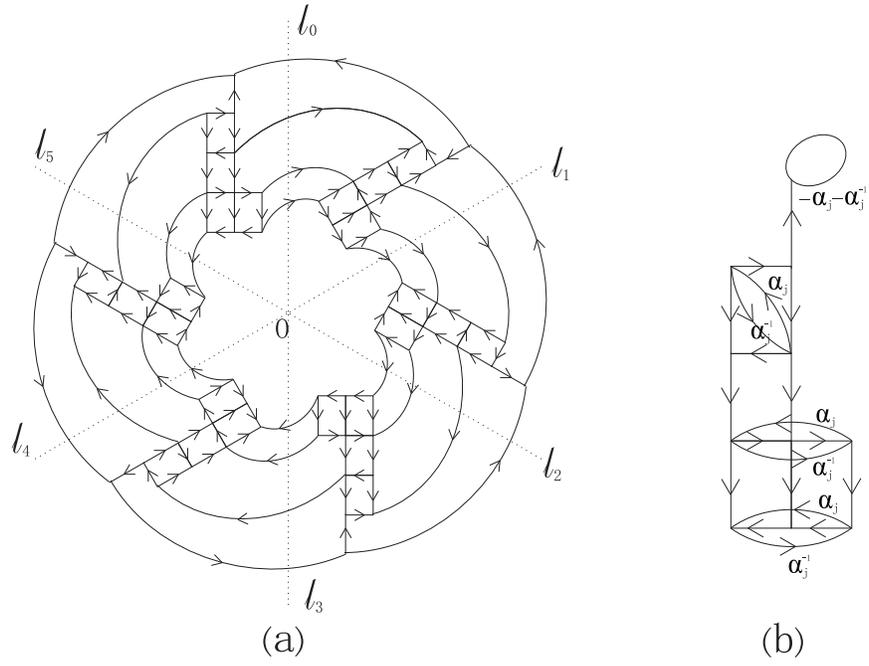}}
  \caption{\ (a)\ A regular orientation $G^e$ of $G$ with local width $w(G)=2$.\ (b)\ The weighted digraph
  $D_j$.}
   %where each arc $(v_s^{(0)},v_{t}^{(0)})$ in $D_j$ which is also be an arc in
%$G_0^e$ is replaced with  two arcs $(v_s^{(0)},v_{t}^{(0)})$ and
%$(v_t^{(0)},v_{s}^{(0)})$ with weights $1$ and $-1$,
%respectively.}
\end{figure}
%%%%%%%%%%%%% Figure 4
%%%%%%%%%%%%%%%%%%%%%%%%%%%%%%%%%%%%%%%%%%
\begin{figure}[htbp]
  \centering
 \scalebox{0.7}{\includegraphics{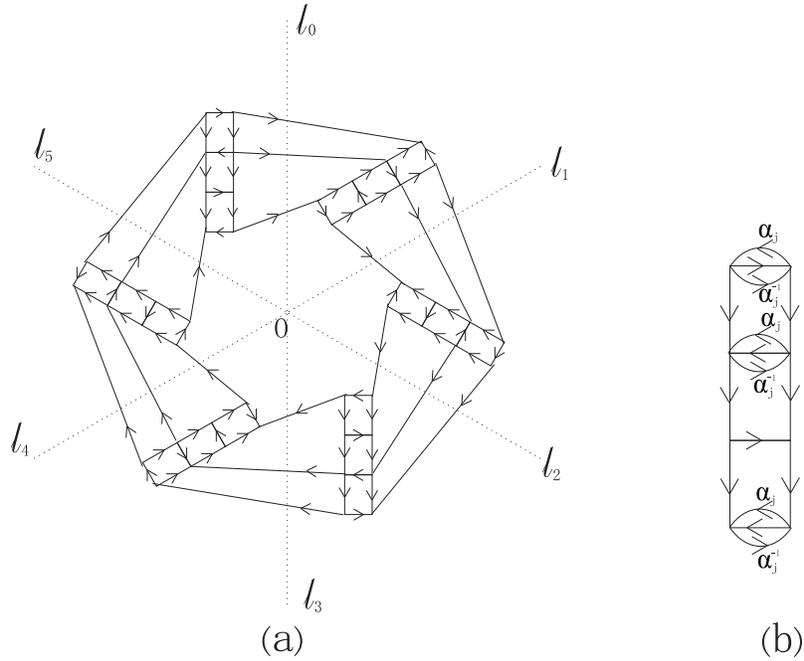}}
  \caption{\ (a)\ A regular orientation $G^e$ of $G$ with the local width $w(G)=1$.
  \ (b)\ The corresponding weighted digraph $D_j$.}
\end{figure}

Hence we have obtained an orientation $G^e$ of $G$. For the graph
$G$ in Figure 1(a) and the orientation $G_0^e$ of $G_0$ in Figure
2(a), by the above definition, the corresponding orientation $G^e$
of $G$ is shown in Figure 3(a). For the underlying graph $G$ of
the digraph illustrated in Figure 4(a), the corresponding
orientation $G^e$ is shown in Figure 4(a).  Note that $G$ is a
bipartite graph. Every boundary face of the $2n$ $G_i$'s in $G^e$
is oddly oriented clockwise. Moreover, it is not difficult to see
that every boundary face in $G^e$ intersected by the rotation axes
$\ell_1,\ell_2,\ldots,\ell_{2n-1}$--except the face containing the
rotation center $O$--is oddly oriented clockwise. For {\it Cases
(i)}, {\it (ii)}, and {\it (iii)}, the number of edges oddly
oriented clockwise in the face containing the rotation center $O$
equals $w(G)\times n+2n-1=1$ (mod $2$) (For the orientation
illustrated in Figure 3(a), the number of edges oddly oriented
clockwise in the face containing the rotation center $O$ equals
$2\times 3+5=11$). For {\it Case (iv)}, the number of edges oddly
oriented clockwise in the face containing the rotation center $O$
equals $w(G)\times n+2n=1$ (mod $2$) (For the orientation shown in
Figure 4(a), the number of edges oddly oriented clockwise in the
face containing the rotation center $O$ equals $1\times 3+6=9$).
Thus the orientation $G^e$ of $G$ satisfies the condition in Lemma
2.3. Consequently, $G^e$ is a Pfaffian orientation of $G$. We call
$G^e$ a regular orientation of $G$.

For $j=0,1,\ldots,n-1$, define $\alpha_j(G)$ as
$$
\alpha_j(G)= \left\{\begin{array}{ll} \cos \frac{2j\pi}{2n}+i\sin
\frac{2j\pi}{2n} &
{\mbox if}\ n=1\ (mod\ 2),\ w(G)=0\ (mod\ 2);\\
\cos \frac{(2j+1)\pi}{2n}+i\sin \frac{(2j+1)\pi}{2n} &
%{\mbox if}\ n=0\ (mod\ 2).
{\mbox otherwise}.
\end{array}
\right. \eqno{(1)}
$$
For convenience, set
$\alpha_j=\alpha_j(G)$ for $0\leq j\leq n-1$.

Now we construct $n$ weighted digraphs $D_0, D_1, \ldots, D_{n-1}$
of order $N/2n$ from $G^e$, which may contain both loops and
multiple arcs, as follows. Note that $G_0^e$ is a digraph
containing neither loops nor multiple arcs and the vertex set of
$G_0^e$ is $V(G_0)=\{v_1^{(0)},v_2^{(0)},\ldots,v_{N/2n}^{(0)}\}$.
Let $\bar {G_0^e}$ be the weighted digraph obtained from $G_0^e$
by replacing each arc $(v_i^{(0)},v_j^{(0)})$ with two arcs
$(v_i^{(0)},v_j^{(0)})$ and $(v_j^{(0)},v_i^{(0)})$ with weights
$1$ and $-1$, respectively. For each $j=0, 1, 2, \ldots, n-1$,
define $D_j$ to be the weighted digraph obtained from
$\bar{G_0^e}$ by the following procedures:

$(1)$\ For each $i=1, 2, \ldots, p$, if $v_{s_i}^{(0)}=
v_{r_i}^{(0)}$, then add a loop at vertex $v_{s_i}^{(0)}$ with
weight $\alpha_j+\alpha_j^{-1}$ if $(v_{r_i}^{(0)},v_{s_i}^{(1)})$
is an arc in $G^e$ and with weight $-\alpha_j-\alpha_j^{-1}$ if
$(v_{s_i}^{(1)},v_{r_i}^{(0)})$ is an arc in $G^e$.

$(2)$\ For each $i=1, 2, \ldots, p$, if $v_{s_i}^{(0)}\neq
v_{r_i}^{(0)}$, then add two arcs $(v_{r_i}^{(0)},v_{s_i}^{(0)})$
and $(v_{s_i}^{(0)},v_{r_i}^{(0)})$ with weights $\alpha_j$ and
$\alpha_j^{-1}$ if $(v_{r_i}^{(0)},v_{s_i}^{(1)})$ is an arc in
$G^e$ and with weights $-\alpha_j$ and $-\alpha_j^{-1}$ if
$(v_{s_i}^{(1)},v_{r_i}^{(0)})$ is an arc in $G^e$.

For the graph $G$ illustrated in Figure 1(a), the corresponding
weighted digraph $D_j$ is shown in Figure 3(b), where each arc
$(v_s^{(0)},v_{t}^{(0)})$ of $D_j$, which is also be an arc in
$G_0^e$, must be regarded as two arcs $(v_s^{(0)},v_{t}^{(0)})$
and $(v_t^{(0)},v_{s}^{(0)})$ with weights $1$ and $-1$,
respectively. For the orientation $G^e$ shown in Figure 4(a), the
corresponding digraph $D_j$ is illustrated in Figure 4(b).

By a suitable labelling of vertices of the regular orientation
$G^e$ of $G$, for {\it Cases (i)-(iii)}, the skew adjacency matrix
$A_1(G^e)$ of $G^e$ has the following form:
$$A_1(G^e)=\left(
\begin{array}{cccccc}
A & R & 0 & \cdots & 0 & R^T\\
-R^T & -A & R & \cdots & 0 & 0 \\
0 & -R^T & A & \cdots & 0 & 0\\
\vdots & \vdots & \ddots & \ddots & \ddots & \vdots\\
0 & 0 & 0 & \cdots & A & R\\
-R & 0 & 0 & \cdots & -R^T & -A
\end{array}
\right)_{2n\times 2n},
$$
and for {\it Case (iv)}, the skew adjacency matrix $A_2(G^e)$ of
$G^e$ has the following form:
$$A_2(G^e)=\left(
\begin{array}{cccccc}
A & R & 0 & \cdots & 0 & -R^T\\
-R^T & -A & R & \cdots & 0 & 0 \\
0 & -R^T & A & \cdots & 0 & 0\\
\vdots & \vdots & \ddots & \ddots & \ddots & \vdots\\
0 & 0 & 0 & \cdots & A & R\\
R & 0 & 0 & \cdots & -R^T & -A
\end{array}
\right)_{2n\times 2n},
$$
where $A$ is the skew adjacency matrix of $G_0^e$, and $R$ is the
adjacent relation in $G_{0,1}^e$ between $G_0$ and $G_1$.

By the definitions of $\alpha_j$'s and $D_j$'s, it is not
difficult to see that the adjacency matrix $A_j$ of weighted
digraph $D_j$ defined above equals exactly
$A+\omega^jR+\omega^{-j}R^T$ if $n=1(\mbox{mod}\ 2)$ and $w(G)$ is
even, and $A+\omega_jR+\omega_j^{-1}R^T$ otherwise, where
$\omega=\cos\frac{2\pi}{2n}+i\sin\frac{2\pi}{2n}$ and
$\omega_j=\cos\frac{(2j+1)\pi}{2n}+i\sin\frac{(2j+1)\pi}{2n}$.

Now we can give a new presentation of the product theorem obtained
independently by Jockusch \cite{Jock94} and by Kuperberg
\cite{Kupe98} as follows.

%%%%%%%%%%%%%%%%%%%%%%%%%%%%%%%%%%%%%%%%%%
%%%%%%%%%%%%%Theorem 3.1
\begin{theorem}[Product Theorem]
Let $G$ be a simple connected plane bipartite graph of order $N$
with $2n$-rotation symmetry and $G_i$'s for $0\leq i\leq 2n-1$ be
the $2n$ graphs defined above. Suppose there are $w(G)$ edges in
$G_i$ lying on the boundary face of $G$ which contains the
rotation center $O$ and there exists no vertex of $G$ lying on the
rotation axes $l_0$, $l_1$, $\ldots$, $l_{2n-1}$ passing $O$. Let
$G^e$ be a regular orientation and $D_j$'s be the $n$ weighted
digraphs defined above. Then the number of perfect matchings of
$G$ can be expressed by
$$M(G)=\prod_{j=0}^{n-1}|\det(A_j)|,$$
where $A_j=(a_{st}^{(j)})$ is the adjacency matrix of $D_j$ each
entry $a_{st}^{(j)}$ of which equals the sum of weights of all
arcs from vertices $v_{s}^{(0)}$ to $v_{t}^{(0)}$ in $D_j$, and
$\alpha_j=\alpha_j(G)$ satisfies $(1)$.
\end{theorem}

%%%%%%%%%%%%%%%%%%%%%%%%%%%%%%%%%%%%%%%%%%
%%%%%%%%%%%%%Corollary
\begin{corollary}
Let $G$ be a bulk plane bipartite lattice of order $N$ with
$2n$-natation symmetry (where $N\rightarrow \infty$). If
$N'=\frac{N}{2n}$ is a constant, then the entropy of $G$
$$
\lim_{N\rightarrow \infty}\frac{2}{N}\log
M(G)=\frac{1}{N'\pi}\int_{0}^{\pi}\log |\det(D)|d_x,
$$
where $D$ is the matrix of order $N'$ obtained from $A_j$ by
replacing each $\alpha_j$ in entries of $A_j$ with $\cos x+i\sin
x$, $i^2=-1$.
\end{corollary}
\begin{proof}
By Theorem 3.1, we have
$$M(G)=\prod_{j=0}^{n-1}|\det(A_j)|,$$
where $A_j$ is the adjacency matrix of $D_j$.% and
%$$\alpha_j(G)=
%\left\{\begin{array}{ll} \cos \frac{j\pi}{n}+i\sin \frac{j\pi}{n}
%&
%{\mbox if}\ n=1\ (mod\ 2),\ w(G)=0\ (mod\ 2);\\
%\cos \frac{(2j+1)\pi}{2n}+i\sin \frac{(2j+1)\pi}{2n} &
%%{\mbox if}\ n=0\ (mod\ 2).
%{\mbox otherwise}.
%\end{array}
%\right. $$
Hence
$$
\lim_{N\rightarrow \infty}\frac{2}{N}\log M(G)=\lim_{n\rightarrow
\infty}\frac{2}{2nN'}\sum_{j=0}^{n-1}\log
|\det(A_j)|=\lim_{n\rightarrow
\infty}\frac{1}{nN'}\sum_{j=0}^{n-1}\log |\det(A_j)|.
$$
By the definition of $A_j$'s, it is not difficult to see that
$$
\lim_{n\rightarrow \infty}\frac{1}{n}\sum_{j=0}^{n-1}\log
|\det(A_j)|=\frac{1}{\pi}\int_{0}^{\pi}\log|\det(D)|d_x
$$which implies the corollary.
\end{proof}

Similarly, we have the following
\begin{corollary}
Let $G$ be a bulk plane bipartite lattice of order $N$ with
$2n$-natation symmetry, where $N'=\frac{N}{2n}\rightarrow \infty$.
Then the entropy of $G$
$$
\lim_{N\rightarrow \infty}\frac{2}{N}\log
M(G)=\frac{1}{\pi}\int_{0}^{\pi}\left[\lim_{N'\rightarrow
\infty}\frac{\log |\det(D)|}{N'}\right]d_x,
$$
where $D$ is the same as in Corollary 3.5.
\end{corollary}

%%%%%%%%%%%%%%%%%%%%%%%%%%%%%%%%%%%%%%%%%%%%%%%%%%%%%%%%%%%%%%%%%%%%%%%%

\section{TWO TYPES OF TILINGS OF CYLINDERS}

Two bulk lattice graphs, denoted by $G_1^*(m,2n)$ and
$G_2^*(m,2n)$, are illustrated in Figure 5(a) and Figure 5(b),
respectively, where $G_1^*(m,2n)$ is a finite subgraph of an
edge-to-edge tilings of the plane with two types of
vertices---8.8.6 and 8.8.4 vertices, and $G_2^*(m,2n)$ is a finite
subgraphs of 8.8.4 tilings in the Euclidean plane which has been
used to describe phase transitions in the layered hydrogen-bonded
$SnCl^2\cdot 2H_2O$ crystal \cite{SN74} in physical systems
\cite{Allen74,NYB89,SN74}. The bulk lattice graph $G_1^*(m,2n)$ is
composed of $2mn$ hexagons whose fundamental part is a hexagon.
Similarly, The bulk lattice graph $G_2^*(m,2n)$ is composed of
$2mn$ quadrangles whose fundamental part is a quadrangle.
Physicists call each of this kind of bulk graphs ``an
$(m,2n)$-bipartite graphs with the free-boundary condition whose
fundamental domain is $G$" (see \cite{LW99}).
%%%%%%%%%%%%%%%%%%%%%%%%%%%%%%%%%%%%%%%%%%
%%%%%%%%%%%%% Figure 5
%%%%%%%%%%%%%%%%%%%%%%%%%%%%%%%%%%%%%%%%%%
\begin{figure}[htbp]
  \centering
 \scalebox{0.7}{\includegraphics{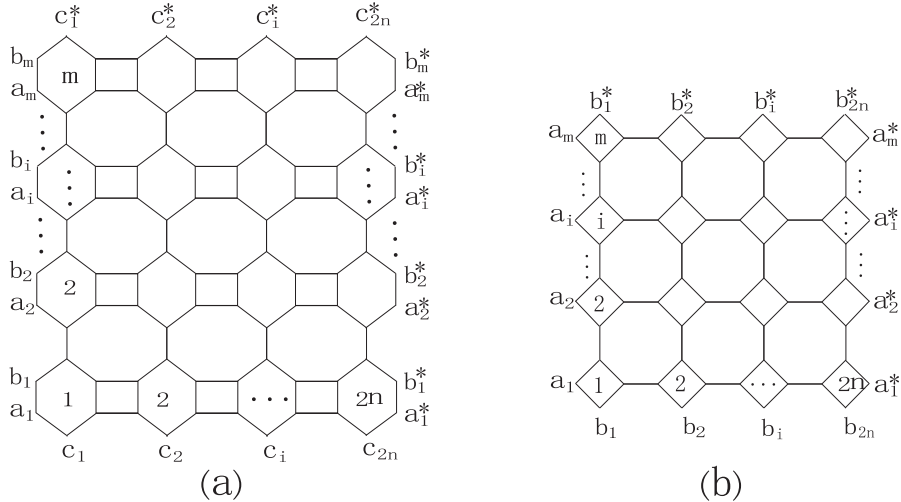}}
  \caption{\ (a)\ The lattice graph $G_1^*(m,2n)$.\ (b)\ The lattice graph $G_2^*(m,2n)$.}
\end{figure}

If we add edges $(a_i,a_i^*), (b_i,b_i^*)$ for $1\leq i\leq m$ and
$(c_i,c_i^*)$ for $1\leq i\leq 2n$ in $G_1^*(m,2n)$, we obtain an
$(m,2n)$-bipartite graph with the doubly-periodic condition on a
torus (see the definition in \cite{KOS06AM}), denoted by
$G_1^t(m,2n)$. Similarly, if we add edges $(a_i,a_i^*)$ for $1\leq
m$ and $(b_i,b_i^*)$ for $1\leq i\leq 2n$ in $G_2^*(m,2n)$, then
an $(m,2n)$-bipartite graph with the doubly-periodic condition on
a torus, denoted by $G_2^t(m,2n)$, is obtained. For some related
work on the plane bipartite graph with the doubly-periodic
condition, see Kenyon, Okounkov, and Sheffield \cite{KOS06AM} and
Cohn, Kenyon, and Propp \cite{CKP01}. Salinas and Nagle
\cite{SN74} showed that the entropy of $G_2^t(m,2n)$, denoted by
$\lim\limits_{n,m\rightarrow \infty}\frac{2}{8mn}\log
[M(G_2^t(m,2n))]$, equals
$$\frac{1}{2\pi}\int_{0}^{\pi/2}\log
\left[\frac{5+\sqrt{25-16\cos^2\theta}}{2}\right]d_{\theta}\approx
0.3770. \eqno{(2)}$$

In this section, as applications of the product theorem in Section
3, we enumerate perfect matchings of two types of subgraphs of
tilings of cylinders---$G_1(m,2n)$ and $G_2(m,2n)$, where
$G_1(m,2n)$ (resp. $G_2(m,2n)$) is obtained from $G_1^*(m,2n)$
(resp. $G_2^*(m,2n)$) by adding extra edges
$(a_i,a_i^*),(b_i,b_i^*)$ for $1\leq i\leq m$ (resp. $(a_i,a_i^*),
(b_j,b_j^*)$ for $1\leq i\leq m, 1\leq j\leq 2n$) between each
pair of opposite vertices of both sides of them. We call each of
$G_1(m,2n)$ and $G_2(m,2n)$  `` an $(m,2n)$-bipartite graph with
the cylinder-boundary condition" (simply cylinder). We also obtain
the exact solutions for the entropies of two types of
corresponding tilings of cylinders. We observe that both the
Archimedean 8.8.4 tilings of the Euclidean plane and the
corresponding tilings of the cylinder have the same entropy which
is approximately 0.3770. Based on the result obtained in
\cite{KOS06AM}, we also show that both $G_1(m,2n)$ and
$G_1^t(m,2n)$ have the same entropy which approximately 0.3344.
These are reasonable conclusions from the physical intuition.

\subsection{THE CYLINDER $G_1(m,2n)$}

Note that the cylinder $G_1(m,2n)$ can be regarded as a
$2n$-rotation symmetric plane bipartite graph. There exists a
regular orientation $G_1(m,2n)^e$ of $G_1(m,2n)$ as stated in
Section 3, which is illustrated in Figure 6(a). For $G_1(m,2n)^e$,
all hexagons in the first column have the same orientation, all
hexagons in the second column have the inverse of the orientation
of hexagons in the first column, and so on.
%%%%%%%%%%%%%%%%%%%%%%%%%%%%%%%%%%%%%%%%%%
%%%%%%%%%%%%% Figure 6
%%%%%%%%%%%%%%%%%%%%%%%%%%%%%%%%%%%%%%%%%%
\begin{figure}[htbp]
  \centering
 \scalebox{0.7}{\includegraphics{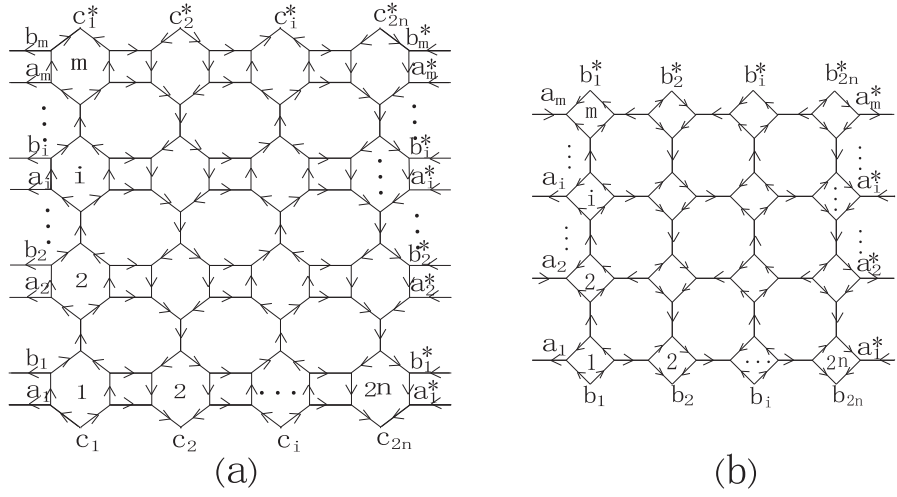}}
  \caption{\ (a)\ The regular orientation $G_1(m,2n)^e$ of $G_1(m,2n)$.
  \ (b)\ The regular orientation $G_2(m,2n)^e$ of $G_2(m,2n)$.}
\end{figure}
%%%%%%%%%%%%%%%%%%%%%%%%%%%%%%%%%%%%%%%%%%%%%%%%%%%%%%%%%%%%%%%%%%%%%%%%
%%%%%%%%%%%%%%%%%%%%%%%%%%%%%%%%%%%%%%%%%%
%%%%%%%%%%%%% Figure 7
%%%%%%%%%%%%%%%%%%%%%%%%%%%%%%%%%%%%%%%%%%
\begin{figure}[htbp]
  \centering
 \scalebox{0.7}{\includegraphics{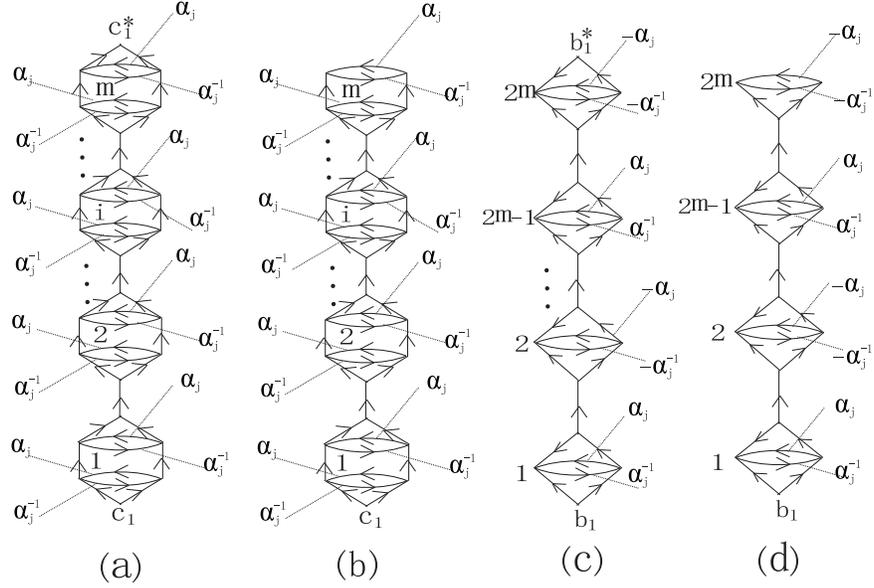}}
  \caption{\ (a)\ The digraphs $D_j$'s corresponding to the regular orientation $G_1(m,2n)^e$ in Figure 6(a).
  \ (b)\ The digraph $D_j'$ obtained from $D_j$ by deleting vertex
  $c_1^*$ shown in Figure 7(a).\ (c)\ The digraphs $D_j$'s corresponding to the regular
orientation $G_2(m,2n)^e$ in Figure 6(b). \ (d)\ The digraph
$D_j'$ obtained from $D_j$ by deleting vertex
  $b_1^*$ shown in Figure 7(c).}
\end{figure}

%%%%%%%%%%%%%%%%%%%%%%%%%%%%%%%%%%%%%%%%%%
%%%%%%%%%%%%%Theorem 4.1
\begin{theorem}
For the cylinder $G_1(m,2n)$, the number of perfect matchings of
$G_1(m,2n)$ can be expressed by
$$M(G_1(m,2n))=\frac{1}{2^n}\prod_{j=0}^{n-1}\frac{1}{\sqrt{4+\beta_j^2}}
\left[\left(\sqrt{4+\beta_j^2}+\beta_j\right)^{2m+1}+\left(\sqrt{4+\beta_j^2}-\beta_j\right)^{2m+1}\right],
\eqno{(3)}$$ and the entropy of $G_1(m,2n)$, {\it i.e., }
$\lim\limits_{m,n\rightarrow \infty}\frac{2}{12mn}\log
M(G_1(m,2n))$,  equals
%$$\lim_{m,n\rightarrow \infty}\frac{2}{12mn}\ln M(G_1(m,2n))=
$$\frac{2}{3\pi}\int_{0}^{\frac{\pi}{2}}\log(\cos x+\sqrt{4+\cos^2
x})d_x\approx 0.3344,$$ where $\beta_j=\cos\frac{j\pi}{n}$ if $n$
is odd and $\beta_j=\cos \frac{(2j+1)\pi}{2n}$ otherwise.
\end{theorem}

\begin{proof}
For the regular orientation $G_1(m,2n)^e$ of the cylinder
$G_1(m,2n)$ shown in Figure 6(a), by the definition of $D_j$ for
$0\leq j\leq n-1$ defined above, the corresponding digraphs
$D_j$'s have the form illustrated in Figure 7(a), where
$\alpha_j$'s satisfy $(1)$, and each arc $(v_i,v_j)$ in $D_j$
whose weight is neither $\alpha_j$ nor $\alpha_j^{-1}$ must be
regarded as two arcs $(v_i,v_j)$ and $(v_j,v_i)$ with weights $1$
and $-1$, respectively. Let $D_j'$ be the digraph obtained from
$D_j$ by deleting vertex $c_1^*$ (see Figure 7(b)).

For $j=0,1,\ldots,n-1$, set
$$L_m(j)=\det(A_j),\ L_m'(j)=\det(A_j'),$$
where $A_j$ (resp. $A_j'$) is the adjacency matrix of the digraph
$D_j$ (resp. $D_j'$). It is not difficult to prove that
$\{L_m(j)\}_{m\geq 0}$ and $\{L_m'(j)\}_{m\geq 0}$ satisfy the
following recurrences:
$$
\left\{\begin{array}{ll}
L_m(j)=(4+4\beta_j^2)L_{m-1}(j)+4\beta_jL_{m-1}'(j) & \mbox{for}\
m\geq 1,\\
L_m'(j)=4\beta_jL_{m-1}(j)+4L_{m-1}'(j) & \mbox{for}\ m\geq 1,\\
L_0(j)=1,\ L_{0}'(j)=0. &
\end{array}
\right.
$$
Hence we have the following:
$$\left\{\begin{array}{ll}
L_m(j)=(8+4\beta_j^2)L_{m-1}(j)-16L_{m-2}(j) & \mbox{for}\ m\geq
2,\\
L_0(j)=1,\ L_1(j)=4+4\beta_j^2,& \end{array}\right.$$ which
implies the following:
$$L_m(j)=\frac{1}{2\sqrt{4+\beta_j^2}}\left[\left(
\sqrt{4+\beta_j^2}+\beta_j\right)^{2m+1}+\left(\sqrt{4+\beta_j^2}-\beta_j\right)^{2m+1}\right].
\eqno{(4)}$$ Hence the equality $(3)$ follows from the product
theorem (Theorem 3.1) and $(4)$. So the entropy of $G_1(m,2n)$
$$\lim\limits_{m,n\rightarrow \infty}\frac{2}{12mn}\log
M(G_1(m,2n))=\lim\limits_{m,n\rightarrow
\infty}\frac{1}{6mn}\times$$$$\left\{-n\log
2-\frac{1}{2}\sum_{j=0}^{n-1}\log(4+\beta_j^2)+\sum_{j=0}^{n-1}
\log\left[\left(\sqrt{4+\beta_j^2}+\beta_j\right)^{2m+1}+
\left(\sqrt{4+\beta_j^2}-\beta_j\right)^{2m+1}\right]\right\}
$$
$$
\begin{array}{lll}
&=&\lim\limits_{m,n\rightarrow
\infty}\frac{1}{6mn}\sum\limits_{j=0}^{n-1}
\log\left[\left(\sqrt{4+\beta_j^2}+\beta_j\right)^{2m+1}+
\left(\sqrt{4+\beta_j^2}-\beta_j\right)^{2m+1}\right]\\
&=&\lim\limits_{m,n\rightarrow
\infty}\frac{1}{3mn}\sum\limits_{j=0}^{[\frac{n-1}{2}]}
\log\left[\left(\sqrt{4+\beta_j^2}+\beta_j\right)^{2m+1}+
\left(\sqrt{4+\beta_j^2}-\beta_j\right)^{2m+1}\right]\\
&=&\lim\limits_{m,n\rightarrow
\infty}\frac{2m+1}{3mn}\sum\limits_{j=0}^{[\frac{n-1}{2}]}\log(\sqrt{4+\beta_j^2}+\beta_j)\\
&=&\frac{2}{3\pi}\int_{0}^{\pi/2}\log(\cos x+\sqrt{4+\cos^2
x})d_x\approx 0.3344
\end{array}
$$
and the theorem thus follows.
\end{proof}

In order to prove the following corollary, we need to introduce a
formula of the entropy for an $(n,n)$-bipartite graphs with the
doubly-period condition obtained by Kenyon, Okounkov, and
Sheffield \cite{KOS06AM}. Let $G$ be a $Z^2$-period bipartite
graph which is embedded in the plane so that translations in the
plane act by color-preserving isomorphisms of $G$---isomorphisms
which map black vertices to black vertices and white to white. Let
$G_n$ be the quotient of $G$ by the action of $nZ^2$. Then $G_n$
is a bipartite graph with the doubly-period condition. Let
$P(z,w)$ be the characteristic polynomial of $G$ (see the
definition in page 1029 in \cite{KOS06AM}). Authors in
\cite{KOS06AM} showed that the entropy of $G_n$
$$\lim_{n\rightarrow \infty}\frac{2}{n^2|G_1|}\log M(G_n)=
\frac{2}{|G_1|(2\pi i)^2}\int_{D}\log|P(z,w)|
\frac{d_z}{z}\frac{d_w}{w}, \eqno{(5)}$$ where $D=\{(z,w)\in C^2:
|z|=|w|=1\}$ and $i^2=-1$.
%%%%%%%%%%%%%%%%%%%%%%%%%%%%%%%%%%%%%%%%%%
%%%%%%%%%%%%%Corollary 4.2
\begin{corollary}
Both $G_1(m,2n)$ and $G_1^t(m,2n)$ have the same entropy, that is,
$$\lim_{m,n\rightarrow \infty}\frac{2}{12mn}\log(M(G_1(m,2n)))=
\lim_{m,n\rightarrow
\infty}\frac{2}{12mn}\log(M(G_1^t(m,2n)))\approx 0.3344.$$
\end{corollary}
%%%%%%%%%%%%%%%%%%%%%%%%%%%%%%%%%%%%%%%%%%%%%%%%%%%%%%%%%%%%%%%%%%%%%%%%
%%%%%%%%%%%%%%%%%%%%%%%%%%%%%%%%%%%%%%%%%%
%%%%%%%%%%%%% Figure 8
%%%%%%%%%%%%%%%%%%%%%%%%%%%%%%%%%%%%%%%%%%
\begin{figure}[htbp]
  \centering
 \scalebox{0.7}{\includegraphics{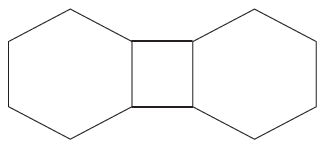}}
  \caption{The fundamental domain of $G_1^t(m,2n)$.}
\end{figure}
\begin{proof}
Note that by the definition in \cite{KOS06AM} the fundamental
domain of $G_1^t(m,2n)$ is composed of two hexagons (see Figure
8). Otherwise, if we use a hexagon as the fundamental domain, then
it does not satisfy the condition ``color-preserving
isomorphisms". It is not difficult to show that the characteristic
polynomial of $G_1^t(m,2n)$
$$P(z,w)=10-4(z+z^{-1})-(w+w^{-1}).$$
Hence, by (5), we have
$$\lim_{m,n\rightarrow\infty}\frac{2}{12mn}\log(M(G_1^t(m,2n)))
$$
$$
\begin{array}{lll}
&=&\frac{2}{12(2\pi)^2}\int_{0}^{2\pi}\int_{0}^{2\pi}\log(10-8\cos
x-2\cos y)d_xd_y\\
&=&\frac{1}{24\pi^2}\int_{0}^{2\pi}\int_{0}^{2\pi}\log(10-8\cos
x-2\cos y)d_xd_y.
\end{array}
$$
Let $F(y)=\int_{0}^{2\pi}\log(10-8\cos x-2\cos
y)d_x=2\int_{0}^{\pi}\log(10-8\cos x-2\cos y)d_x$. Then
$$F'(y)=4\sin y\int_{0}^{\pi}\frac{d_x}{10-8\cos x-2\cos y}=\frac{2\pi\sin y}{\sqrt{\cos^2y-10\cos y+9}}.$$
Hence we have
$$F(y)=2\pi\log(5-\cos y+\sqrt{\cos^2y-10\cos y+9})$$
implying
$$\frac{1}{24\pi^2}\int_{0}^{2\pi}\int_{0}^{2\pi}\log(10-8\cos
x-2\cos y)d_xd_y
$$$$
\begin{array}{lll}
&=&\frac{1}{12\pi}\int_{0}^{2\pi}\log(5-\cos
x+\sqrt{\cos^2x-10\cos x+9})d_x\\
&=&\frac{1}{6\pi}\int_{0}^{\pi}\log(5-\cos x+\sqrt{\cos^2x-10\cos
x+9})d_x
\end{array}$$
 So, by Theorem 4.1, it suffices to prove
the following:
$$4\int_{0}^{\frac{\pi}{2}}\log(\cos x+\sqrt{4+\cos^2
x})d_x=\int_{0}^{\pi}\log(5-\cos x+\sqrt{\cos^2x-10\cos
x+9})d_x.$$ Note that
$$4\int_{0}^{\frac{\pi}{2}}\log(\cos x+\sqrt{4+\cos^2
x})d_x=2\int_{0}^{\frac{\pi}{2}}\ln(2\cos^2x+4+\sqrt{16\cos^2
x+4\cos^4x})d_x$$
$$=\int_{0}^{\pi}\log(2\cos^2\frac{x}{2}+4\sqrt{16\cos^2\frac{x}{2}+4\cos^4\frac{x}{2}})d_x
=\int_{0}^{\pi}\log(5+\cos x+\sqrt{\cos^2x+10\cos x+9})d_x$$
$$=\int_{0}^{\pi}\log(5-\cos x+\sqrt{\cos^2x-10\cos x+9})d_x.$$
The corollary has thus been proved.
\end{proof}
%%%%%%%%%%%%%%%%%%%%%%%%%%%%%%%%%%%%%%%%%%%%%%%%%%%%%%%%%%%%%%%%%%%%%%%%

\subsection{THE CYLINDER $G_2(m,2n)$}

Note that the cylinder $G_1(m,2n)$ can be embedded into the plane
such that it is a $2n$-rotation symmetric plane bipartite graph.
There exists a regular orientation $G_2(m,2n)^e$ of $G_2(m,2n)$ as
stated in Section 3, which is illustrated in Figure 6(b).
%\newpage
%%%%%%%%%%%%%%%%%%%%%%%%%%%%%%%%%%%%%%%%%%
%%%%%%%%%%%%%Theorem 4.1
\begin{theorem}
For the cylinder $G_2(m,2n)$, the number of perfect matchings of
$G_2(m,2n)$ can be expressed by
$$M(G_2(m,2n))=$$$$
\prod\limits_{j=0}^{n-1}\left[\frac{\sqrt{9+16\beta_j^2}+3}{2\sqrt{9+16\beta_j^2}}
\left(\frac{5+\sqrt{9+16\beta_j^2}}{2}\right)^{m}\right.+
\left.\frac{\sqrt{9+16\beta_j^2}-3}{2\sqrt{9+16\beta_j^2}}
\left(\frac{5-\sqrt{9+16\beta_j^2}}{2}\right)^{m}\right],\eqno{(6)}$$
and the entropy, {\it i.e.}, $\lim\limits_{m,n\rightarrow
\infty}\frac{2}{8mn}\log M(G_2(m,2n))$, equals
$$\frac{1}{2\pi}\int_{0}^{\pi/2}\log
\left[\frac{5+\sqrt{25-16\cos^2\theta}}{2}\right]d_{\theta}\approx
0.3770,$$ where $\beta_j=\cos\frac{j\pi}{n}$ if $n$ is odd and
$\beta_j=\cos \frac{(2j+1)\pi}{2n}$ otherwise.
\end{theorem}

\begin{proof}
We can prove easily the statement in the theorem on the entropy
from $(5)$. Hence it suffices to prove that $(5)$ holds. For the
regular orientation $G_2(m,2n)^e$ of the cylinder $G_2(m,2n)$
shown in Figure 6(b), by the definition of $D_j$ for $0\leq j\leq
n-1$ defined above, the corresponding digraphs $D_j$'s have the
form illustrated in Figure 7(c), where $\alpha_j$'s satisfy $(1)$,
and each arc $(v_i,v_j)$ in $D_j$ whose weight is neither
$\alpha_j$ (or $-\alpha_j$) nor $\alpha_j^{-1}$ (or
$-\alpha_j^{-1}$) must be regarded as two arcs $(v_i,v_j)$ and
$(v_j,v_i)$ with weights $1$ and $-1$, respectively. Let $D_j'$ be
the digraph obtained from $D_j$ by deleting vertex $b_1^*$ (see
Figure 7(d)).

For $j=0,1,\ldots,n-1$, set
$$P_m(j)=\det(A_j),\ P_m'(j)=\det(A_j'),$$
where $A_j$ (resp. $A_j'$) is the adjacency matrix of the digraph
$D_j$ (resp. $D_j'$) illustrated in Figure 6(c) (resp. Figure
6(d)). It is not difficult to prove that $\{P_m(j)\}_{m\geq 0}$
and $\{P_m'(j)\}_{m\geq 0}$ satisfy the following recurrences:
$$
\left\{\begin{array}{ll}
P_{2m+1}(j)=4P_{2m}(j)+2\beta_jP_{2m}'(j),\
P_{2m+1}'(j)=-2\beta_jP_{2m}(j)-P_{2m}'(j) &  m\geq 0,\\
P_{2m}(j)=4P_{2m-1}(j)-2\beta_jP_{2m-1}'(j),\
P_{2m}'(j)=2\beta_jP_{2m-1}(j)-P_{2m-1}'(j) &  m\geq 1,\\
P_0(j)=1,\ P_{0}'(j)=0,\ P_1(j)=4,\ P_1'(j)=-2\beta_j. &
\end{array}
\right.
$$
Hence from the recurrences above we have the following:
$$\left(\begin{array}{c}P_{2m+1}(j)\\
P_{2m+1}'(j)\end{array}\right)=\left(\begin{array}{cc}16+4\beta_j^2 & -10\beta_j\\
-10\beta_j & 4\beta_j^2+1\end{array}\right)\left(\begin{array}{c}P_{2m-1}(j)\\
P_{2m-1}'(j)\end{array}\right) \eqno{(7)}$$ and
$$\left(\begin{array}{c}P_{2m}(j)\\
P_{2m}'(j)\end{array}\right)=\left(\begin{array}{cc}16+4\beta_j^2 & 10\beta_j\\
10\beta_j & 4\beta_j^2+1\end{array}\right)\left(\begin{array}{c}P_{2(m-1)}(j)\\
P_{2(m-1)}'(j)\end{array}\right).\eqno{(8)}$$ Let
$a_m=P_{2m+1}(j)$ and $b_m=P_{2m}(j)$ for $j\geq 0$. Hence we have
$$\left\{\begin{array}{ll}
a_m=(18\beta_j^2+17)a_{m-1}-16(1-\beta_j^2)^2a_{m-2} & \mbox{for}\
m\geq
2,\\
b_m=(18\beta_j^2+17)b_{m-1}-16(1-\beta_j^2)^2b_{m-2} & \mbox{for}\
m\geq 2,\\
a_0=4,a_1=64+36\beta_j^2,b_0=1,b_1=16+4\beta_j^2.&
\end{array}\right. \eqno{(9)}$$
By solving the recurrences in $(9)$, then we have

$$
a_m=\prod\limits_{j=0}^{n-1}\left[\frac{\sqrt{9+16\beta_j^2}+3}{2\sqrt{9+16\beta_j^2}}
\left(\frac{5+\sqrt{9+16\beta_j^2}}{2}\right)^{2m+1}\right.+
\left.\frac{\sqrt{9+16\beta_j^2}-3}{2\sqrt{9+16\beta_j^2}}
\left(\frac{5-\sqrt{9+16\beta_j^2}}{2}\right)^{2m+1}\right],$$ and
$$
b_m=\prod\limits_{j=0}^{n-1}\left[\frac{\sqrt{9+16\beta_j^2}+3}{2\sqrt{9+16\beta_j^2}}
\left(\frac{5+\sqrt{9+16\beta_j^2}}{2}\right)^{2m}\right.+
\left.\frac{\sqrt{9+16\beta_j^2}-3}{2\sqrt{9+16\beta_j^2}}
\left(\frac{5-\sqrt{9+16\beta_j^2}}{2}\right)^{2m}\right],
$$
Hence $(6)$ has been proved and the theorem thus follows.
\end{proof}

From $(2)$ and Theorem 4.2, both the 8.8.4 tilings of the
Euclidean plane and the corresponding tilings of the cylinder have
the same entropy which is approximately 0.3770.

%%%%%%%%%%%%%%%%%%%%%%%%%%%%%%%%%%%%%%%%%%%%%%%%%%%%%%%%%%%%%%%%%%%%%%%%

\section{CONCLUDING REMARKS}
In the classical work on the dimer problem for plane quadratic
lattices, Kasteleyn \cite{Kast61} proved that the $m\times n$
quadratic lattice with the free-boundary condition, the $m\times
n$ quadratic lattice with the doubly-period condition ({\it i.e.},
the $m\times n$ torus), and the $m\times n$ quadratic lattice with
the cylinder-boundary condition ({\it i.e.}, the $m\times n$
cylinder) have the same entropy ($=2c/\pi\approx 0.5831$, where
$c$ is Catalan's constant). In this paper we have investigated the
problem on enumeration of perfect matchings of plane bipartite
graphs with $2n$-rotation symmetry. We compute the entropy of a
bulk plane bipartite lattice with $2n$-notation symmetry. We
obtain the explicit expressions for the numbers of perfect
matchings and entropies for two types of tilings of (the surface
of) cylinders and we showed that both the bipartite graph
$G_1(m,2n)$ with the cylinder-period condition and the bipartite
graph $G_1^t(m,2n)$ with the doubly-period condition have the same
entropy and that both the bipartite graph $G_2(m,2n)$ with the
cylinder-period condition and the bipartite graph $G_2^t(m,2n)$
with the doubly-period condition have the same entropy. A natural
problem is how to enumerate perfect matchings of the plane
bipartite graph with $(2n+1)$-rotation symmetry. A more general
problem is to enumerate perfect matchings of the plane graph (not
necessary to be bipartite) with $n$-rotation symmetry.
%\newpage
%%%%%%%%%%%%%%%%%%%%%%%%%%%%%%%%%%%%%%%%%%%
%%%%%%%%%%%%%Conjecture 5.1
%\begin{conjecture}
%Let $G(m,n), G^t(m,n)$, and $G^*(m,n)$ be the three plane
%bipartite graphs with the cylinder-period condition, the
%doubly-period condition, and the free-boundary condition
%respectively, where their fundamental domain is a plane bipartite
%graph $G$ with perfect matchings. Then $G(m,n), G^t(m,n)$, and
%$G^*(m,n)$ have the same entropy. That is,
%$$\lim_{m,n\rightarrow \infty}\frac{2}{mn|G|}\log M(G(m,n))=$$$$
%\lim_{m,n\rightarrow \infty}\frac{2}{mn|G|}\log M(G^t(m,n))=
%\lim_{m,n\rightarrow \infty}\frac{2}{mn|G|}\log M(G^*(m,n)),$$
%where $|G|$ is the number of vertices of $G$.
%\end{conjecture}
%\newpage
%%%%%%%%%%%%%%%%%%%%%%%%%%%%%%%%%%%%%%%%%%%
%%%%%%%%%%%%%%% Acknowledgements
%%%%%%%%%%%%%%%%%%%%%%%%%%%%%%%%%%%%%%%%%%%
\vskip0.5cm \noindent {\bf Acknowledgements} \vskip0.1in
\par Thanks to Professors G. Kuperberg and C. Krattenthaler for kindly telling
us references \cite{Jock94} and \cite{Kupe98}. We would like to
thank Professor Z. Chen for providing many very helpful
suggestions for revising this paper. The third author Fuji Zhang
thanks Institute of Mathematics, Academia Sinica for its financial
support and hospitality.

\bibliographystyle{amsplain}

\end{document}